# POINCARÉ INEQUALITIES AND QUASICONFORMAL STRUCTURE ON THE BOUNDARY OF SOME HYPERBOLIC BUILDINGS

MARC BOURDON AND HERVÉ PAJOT

ABSTRACT. In this paper we shall show that the boundary $\partial I_{p,q}$ of the hyperbolic building $I_{p,q}$ considered in [1] admits Poincaré type inequalities. Then by using Heinonen-Koskela's work [7], we shall prove Loewner capacity estimates for some families of curves of $\partial I_{p,q}$ and the fact that every quasiconformal homeomorphism $f : \partial I_{p,q} \longrightarrow \partial I_{p,q}$ is quasisymetric. Therefore by these results, the answers to questions 19 and 20 of Heinonen and Semmes [8] are NO.

## 1. INTRODUCTION

In recent work Heinonen and Koskela [7] showed that in metric spaces in which the modulus of the family of curves joining two continua is controlled, quasiconformal homeomorphisms are quasisymetric. They characterized such spaces (called Loewner spaces) by the existence of Poincaré type inequalities. For instance, $\mathbb{R}^n$ ($n \geq 2$), Carnot groups, and so the boundary of any non compact symmetric space of rank 1 (and dimension at least 3) are Loewner spaces. In this paper we shall show that the boundary of some hyperbolic buildings belongs also to this class of spaces.

For every $p \geq 5$ and every $q \geq 3$, we denote by $I_{p,q}$ the Tits building whose apartments are hyperbolic planes with curvature $-1$, whose chambers are regular hyperbolic $p$-gons with angle $\frac{\pi}{2}$ and whose link of each vertex is the complete bipartite graph with $q + q$ vertices (see [2] and [11] for complete treatment about buildings, see [1] and section 2 for more details on hyperbolic buildings).

The hyperbolic building $I_{p,q}$ is a hyperbolic space in the sense of Gromov (see [3]). More precisely it is a CAT($-1$)-space (which means that its triangles are thiner than those of the hyperbolic space with curvature -1) and therefore has a boundary at infinity $\partial I_{p,q}$ which is homeomorphic to the Menger's Universal Curve (recall that this set is a continuum obtained as an analogue of the classical Cantor middle-third set by punching holes out of a cube in some regular manner). We can equip the boundary of $I_{p,q}$ with a natural metric $\delta_{p,q}$ which has the following properties:
- the metric space $\partial I_{p,q}$ is geodesic (which means that every couple of points can be joined by a curve whose length is equal to the distance between the two points);
- its Hausdorff dimension $Q_{p,q}$ is equal to Pansu conformal dimension of $\partial I_{p,q}$;

1991 *Mathematics Subject Classification*. Primary: 30C65, 51E24.

Parts of this work were done during a stay of the second author at MSRI. Research at MSRI is supported in part by NSF grant DMS-9022140.

- its Hausdorff measure $\mu_{p,q}$ is Ahlfors-regular with dimension $Q_{p,q}$, which means that there exists $C > 0$ such that, for every $\xi \in \partial I_{p,q}$ and every $r \in (0, \operatorname{diam} \partial I_{p,q})$,

$$(1) \qquad C^{-1} r^{Q_{p,q}} \leq \mu_{p,q}(\partial I_{p,q} \cap B(\xi, r)) \leq C r^{Q_{p,q}}.$$

Throughout the paper, $B(\xi, r)$ denotes the open ball in $\partial I_{p,q}$ with center $\xi$ and radius $r$ for the metric $\delta_{p,q}$.

*Remark.* The number $Q_{p,q} = 1 + \dfrac{\log(q-1)}{\operatorname{Argch}(\frac{p-2}{2})}$ (see [1], Theorem 1.1) is strictly bigger than 1 and is not in general an integer.

We shall prove in this paper the following results:

**Theorem 1.** *Let $I_{p,q}$ be a hyperbolic building as defined above.*
*Then its boundary $\partial I_{p,q}$ admits weak $(1, \alpha)$-Poincaré inequalities for every $\alpha \geq 1$, which means that there exist constants $C_0 > 0$ and $C_\alpha > 0$ such that*

$$(2) \qquad \fint_B |u - u_B| d\mu_{p,q} \leq C_\alpha \operatorname{diam} B \left( \fint_{C_0 B} \rho^\alpha d\mu_{p,q} \right)^{\frac{1}{\alpha}}$$

*whenever*
- *$B$ is an open ball in $\partial I_{p,q}$;*
- *$u : \partial I_{p,q} \longrightarrow \mathbb{R}^+$ is a bounded continuous function in the ball $C_0 B$ (which is the ball with the same center as $B$ but whose radius is equal to $C_0$ times the radius of $B$);*
- *$\rho : \partial I_{p,q} \longrightarrow \mathbb{R}^+$ is a very weak gradient of $u$ in $C_0 B$, which means that for any two points $x$, $y$ in $C_0 B$,*

$$(3) \qquad |u(x) - u(y)| \leq \int_\gamma \rho(s) ds$$

*for every rectifiable curve $\gamma$ joining $x$ and $y$ in $C_0 B$.*
*We denote by $u_B$ the average of $u$ in $B$ : $u_B = \dfrac{1}{\mu_{p,q}(B)} \int_B u d\mu_{p,q} = \fint_B u d\mu_{p,q}$.*

*Remarks.* (i) The weak $(1,1)$-Poincaré inequality is the strongest Poincaré inequality as it implies the others by Hölder inequalities. Therefore we shall prove inequality (2) only for $\alpha = 1$. It should be noted that the Poincaré inequalities (2) are similar to the usual inequalities in $\mathbb{R}^n$ (see e.g. [5]).
(ii) The definition of the very weak gradient is not so surprising. For instance, in $\mathbb{R}^n$, $|\nabla u|$ is a very weak gradient of the smooth function $u$. It should also be mentioned that for every function $u$ in $\partial I_{p,q}$, a very weak gradient always exists (take $\rho = \infty$) and it is not unique.

**Theorem 2.** *Let $I_{p,q}$ be a hyperbolic building as defined above. Then,*
*a) $\partial I_{p,q}$ is a Loewner space ;*
*b) every quasiconformal homeomorphism $f : \partial I_{p,q} \longrightarrow \partial I_{p,q}$ is quasisymetric.*

This theorem is an easy consequence of the results contained in [7] and of Theorem 1. We now provide some of the definitions needed.

The metric space $\partial I_{p,q}$ is a Loewner space if the modulus of the family of curves joining two continua is controlled. More precisely, for every $t > 0$,

$$(4) \qquad \lambda(t) = \inf\{\mathrm{mod}(E, F, \partial I_{p,q}) : \Delta(E, F) \leq t\} > 0$$

whenever $E$ and $F$ are non degenerate continua in $\partial I_{p,q}$ (recall that a continuum is a compact connected set) ,

$$\Delta(E, F) = \frac{\mathrm{dist}(E, F)}{\min(\mathrm{diam}\ E,\ \mathrm{diam} F)} \quad \text{and}$$

$$\mathrm{mod}(E, F, \partial I_{p,q}) = \inf \int_{\partial I_{p,q}} \rho^{Q_{p,q}} d\mu_{p,q},$$

where the infimum is taken over all the measurable positive functions $\rho$ such that $\int_\gamma \rho(s)ds \geq 1$ for every rectifiable curve $\gamma$ joining $E$ to $F$.

Roughly speaking, in a Loewner space, there are a lot of nice curves joining two continua. Since the metric measure space $\partial I_{p,q}$ is geodesic, proper (which means that closed balls are compact) and Ahlfors regular, the Theorem 2 follows from Theorem 1 and Theorem 5.7 in [7] (which gives a characterization of Loewner spaces in terms of weak Poincaré inequalities).

Let $f : \partial I_{p,q} \to \partial I_{p,q}$ be a homeomorphism.
For every $\xi \in \partial I_{p,q}$, every $r > 0$, set

$$H_f(\xi, r) = \frac{\sup\{|f(\xi) - f(\eta)|; |\xi - \eta| \leq r\}}{\inf\{|f(\xi) - f(\eta)|; |\xi - \eta| \geq r\}}.$$

We say that $f$ is quasiconformal (QC) if there exists $H < \infty$ such that, for every $\xi \in \partial I_{p,q}$,

$$\limsup_{r \to 0} H_f(\xi, r) \leq H.$$

We say that $f$ is quasisymmetric (QS) if there exists $H < \infty$ such that, for every $\xi \in \partial I_{p,q}$, every $r > 0$, $H_f(\xi, r) \leq H$.

Quasiconformal (respectively quasisymmetric) homeomorphisms distort the shape of infinitesimal balls (respectively of every ball) by a uniformly bounded amount. For instance, a quasiconformal (respectively quasisymmetric) homeomorphism of $\mathbb{C}$ transforms a ball of small radius (respectively every ball) in an "ellipse" whose eccentricity is bounded. Theorem 2 part b) is an easy application of Theorem 4.9 of [7] and of the fact that $\partial I_{p,q}$ is a Loewner space. In fact, the equivalence between quasiconformal structures and quasisymmetric structures was proved by Gehring [4] in $\mathbb{R}^2$, by Gehring-Väisälä in $\mathbb{R}^n$, $n \geq 2$ (see [13]), and by Heinonen-Koskela [6] in Carnot groups. The later authors showed (see [7]) that the crucial point to obtain this equivalence was the Loewner capacity estimates (4).

Väisälä [14] proved that every quasisymmetric homeomorphism of a metric space is quasimöbius. Moreover, Paulin [10] showed that every quasimöbius homeomorphism of the boundary of a Gromov quasi-homogeneous hyperbolic space is induced by a quasiisometry of the space. Hence, from this and Theorem 2 we deduce the following result (see [6] for similar theorems for Carnot groups).

**Corollary 3.** *Let $I_{p,q}$ be a hyperbolic building as above and let $f : \partial I_{p,q} \to \partial I_{p,q}$ be a quasiconformal homeomorphism. Then, $f$ is induced by a quasiisometry $F : I_{p,q} \to I_{p,q}$.*

A mapping $F : I_{p,q} \to I_{p,q}$ is a quasiisometry if there exists $\lambda > 1$, $k \geq 0$ such that
$$\lambda^{-1}|x - x'| - k \leq |F(x) - F(x')| \leq \lambda|x - x'| + k,$$
whenever $x$ and $x'$ are in $I_{p,q}$.
*Remark.* The converse of Corollary 3 is true for general hyperbolic spaces (see [3]).

It seems that the boundaries of hyperbolic buildings $I_{p,q}$ are the first examples of regular metric spaces admitting weak $(1,1)$-Poincaré inequalities and of Loewner spaces with Hausdorff dimension $Q$ not an integer. Thus the answer to question 19 of [8] ("If $X$ is an $Q$-Ahlfors regular space that admits a weak $(1,1)$-Poincaré inequality, is then $Q$ an integer ?") and the answer to question 20 ("If $X$ is an Ahlfors regular Loewner space for some $Q > 1$, is $Q$ then an integer ?") are NO.
Our arguments to prove (2) are quite standard (see for instance the discussion in [12]): we shall prove an inequality equivalent to (2). Namely, whenever $\xi, \eta \in \partial I_{p,q}$, we control the variation $|u(\xi) - u(\eta)|$ by some maximal functions of the very weak gradient of $u$. The key point is then to find a "good" family of curves joining two points $\xi$ and $\eta$ in $\partial I_{p,q}$.
The paper is organized as follows. In section 2, we introduce the hyperbolic buildings and we give some of their useful properties. In section 3, we shall prove theorem 1.

## 2. The metric measure space $\partial I_{p,q}$

The hyperbolic building $I_{p,q}$ has the following property. It is the unique simply connected cell 2-complex such that:
- its 2-cells are regular hyperbolic $p$-gons with angles $\frac{\pi}{2}$;
- two of its 2-cells share at most one edge or one vertex;
- the link of each vertex is the complete bipartite graph with $q + q$ vertices.
Recall that the link of the vertex $x$ of the complex $X$ is the graph $L_x$ such that:
- its vertices are the edges of $X$ containing $x$;
- two vertices $i$ and $j$ of $L_x$ are related by an edge if a face of $X$ contains both the edges represented by $i$ and $j$ (see figure 1).
An apartment $A$ of $I_{p,q}$ is a copy of the hyperbolic plane of curvature $-1$. It is equipped with a natural pavage by its chambers. A wall of $A$ is a hyperbolic geodesic contained in the 1-skeleton of the pavage.

FIGURE 1. The building $I_{6,3}$ and its link $L_{3,3}$

A geodesic ray $\mathcal{R}$ of $I_{p,q}$ is defined by $\mathcal{R} = s([0, \infty))$ where $s : \mathbb{R} \to I_{p,q}$ is an isometry. We say that two geodesic rays are asymptotic if their Hausdorff distance is finite. This relation on the set of geodesic rays is an equivalence relation and the set of equivalence classes coincides with the boundary of $I_{p,q}$ (see [3] Proposition 3.2).

Let $I$ be a hyperbolic building as above. We omit henceforth the letters $p$ and $q$ and we shall follow in this section the terminology of [1].

We choose now a basepoint $x$ in the space $I$ and let $\delta_x$ be the metric on the boundary $\partial I$ of $I$ related to $x$ and defined in [1] (Lemma 3.1.4). This metric has the following fundamental property:

There exists a positive constant $C$ so that

$$(5) \qquad C^{-1} a^{-\{\zeta|\zeta'\}_x} \leq \delta_x(\zeta, \zeta') \leq C a^{-\{\zeta|\zeta'\}_x},$$

whenever $\zeta, \zeta'$ are in $\partial I$. Here $\{.|.\}_x$ denotes the combinatorial Gromov product defined in [1] 2.4.C and $a = \exp\left(\dfrac{\log(q-1)}{Q-1}\right)$. In particular, the combinatorial Gromov product $\{\xi|\xi'\}_x$ is the number of walls of $A$ which intersect both $[x, \xi)$ and $[x, \xi')$ whenever $x, \xi$ and $\xi'$ are in the same apartment.

*Remark.* the last inequality is classic for visual metrics and Gromov product in hyperbolic spaces (see [3]).

Thus, fix two points $\xi$ and $\eta$ in $\partial I$. The main goal of this section is to find a nice family of curves $(\Gamma_t)$ joining $\xi$ and $\eta$ in $10B_{\xi,\eta}$ (where $B_{\xi,\eta}$ is a ball in $\partial I$ which contains $\xi$ and $\eta$ and whose diameter is $\delta_x(\xi, \eta)$) such that we can control the size of the set of parameters $t$.

We begin with some notations and remarks. Let $A$ be an apartment of $I$ whose boundary contains $\xi$ and $\eta$, $c$ be a chamber of $A$ which contains $x$, $K$ be the compact subgroup of $\text{Isom}(I)$ which fixes $c$ point by point and $H$ be the subgroup of $K$ which fixes $\xi$ and $\eta$. We can always suppose that $x \in A$, since we have for each $y$ in a neighborhood of $A$, $\delta_y \sim a^{-|y-z|}\delta_z$ where $z$ is the "orthogonal projection" of $y$ on $A$ (see [1] 2.4.5).

Note that $K$ and $H$ act isometrically on $(\partial I, \delta_x)$.

We denote by $dk$ and $dh$ the Haar probabilities of $K$ and $H$ and by $da$ the Hausdorff

measure of $(\partial A, \delta_{x\ |\partial A})$. Consider the continuous and surjective projection $\Pi$

(6) $$\Pi\ :\ K \times \partial A \longrightarrow \partial I$$
(7) $$(k, \zeta) \longrightarrow \Pi_\zeta(k) = k\zeta.$$

Thus, $\mu = \Pi_*(dk \times da)$ (namely the image measure of $dk \times da$ by $\Pi$) is $Q$-regular on $\partial I$ where $Q$ is the Hausdorff dimension of $\partial I$ (see [9] for the definition of Hausdorff dimension and Hausdorff measures).

Finally, we denote by $[\xi\eta]$ the segment of $\partial A$ joining $\xi$ and $\eta$ of smallest length and $\mathcal{C}$ the "cone" defined by

$$\mathcal{C} = \Pi(H \times [\xi\eta]) = \bigcup_{h \in H} h[\xi\eta].$$

The restriction of the metric $\delta_x$ to $[\xi, \eta]$ is geodesic; therefore $[\xi, \eta]$ is isometrically equivalent to the real interval $[0, l]$ where $l = \delta_x(\xi, \eta)$. We denote by $t$ the point of $[\xi, \eta]$ whose distance from $\xi$ is $t$. We identify the measure $da$ with $dt$. For every $t \in [0, l]$, let $\gamma_t$ be the probability $\Pi_{t*}(dk)$ on $Kt$ $(= \{kt, k \in K\}$, namely the fiber over $t)$.

Note that $(h[\xi, \eta])_{h \in H}$ is a family of curves joining $\xi$ and $\eta$ of length $\delta_x(\xi, \eta)$ and the following lemma will give an estimate on the size of the "set of parameters" $H$ on each fiber.

**Lemma 4.** *a) For every $t \in (0, l)$, the probability $\Pi_{t*}(dh)$ on $Ht$ is*

(8) $$\Pi_{t*}(dh) = \frac{1}{\gamma_t(Ht)} \gamma_{t\ |Ht}.$$

*b) There exists a positive constant $C > 0$ (which depends only on $I$) such that*

(9) $$C^{-1}(\phi_l(t))^{Q-1} \leq \gamma_t(Ht) \leq C(\phi_l(t))^{Q-1}$$

*whenever $t \in (0, l)$, and where $\phi_l(t) = \inf(t, l - t)$.*

*Proof.* Note that b) implies a). Indeed, by b), $\gamma_t(Ht) \neq 0$ for every $t \in (0, l)$, therefore the probability $\frac{1}{\gamma_t(Ht)} \gamma_{t\ |Ht}$ is well defined. Moreover, the probabilities $\Pi_{t*}$ and $\frac{1}{\gamma_t(Ht)} \gamma_{t\ |Ht}$ are $H$-invariant on $Ht$, since $Ht$ is homogeneous under $H$, thus they are equal (by unicity of the invariant measure, see for instance [9], Theorem 3.1).

We begin now with the proof of b).
We shall start by defining geometrically $Kt$ and $Ht$. For this, we consider first the geodesic rays $[x\xi)$, $[x\eta)$ and $[xt)$ of $A$ (recall that the geodesic ray $[x, \xi)$ is the image of the interval $[0, \infty)$ by an isometry $s$ such that $s(0) = x$ and $s(\infty) = \xi$). Even if it means moving $x$ in the chamber $c$, we can suppose that $[xt)$ does not pass through any vertex of the pavage of $A$ by its chambers.
Set $T = K[xt)$. Then, $T$ is a tree whose root is $x$ and whose geodesic rays from $x$ are the images of the ray $[x, t)$ by the elements of $K$. Its vertices are the points of intersection of $T$ with the walls of $I$. Since $[xt)$ does not intersect any vertex of the pavage of $A$, the valence of each vertex is equal to $q$. We can identify the boundary of $T$ with $Kt$ and the probability $\gamma_t$ with the standard probability of a rooted tree of valence $q$.

We consider now the subset $U = H[xt)$ of $T$. To describe it, we associate each vertex of the ray $[xt)$ of $T$ with an integer. Since each vertex of $T$ is the image of an unique vertex of $[xt)$ by an unique element of $K$, the vertices of $T$ can be equipped by a $K$-invariant labeling. We say that a number $n \in \mathbb{N}$ is good if the wall of $A$ which intersects $[xt)$ at the vertex "$n$" intersects neither $[x\xi)$ nor $[x\eta)$. Thus, $U$ is a subtree of $T$ whose vertices which have good numbers are of valence $q$ and whose vertices which carry bad numbers are of valence 2. To see this, remark that the apartments which contain the chamber $c$ and whose boundary contains $\xi$ and $\eta$ are obtained in this following way: we "bend" the apartment $A$ along a wall which intercepts neither $[x, \xi)$ nor $[x, \eta)$ whereas the demi-apartment containing $c$ remains fixed so that we obtain a new apartment $A'$ containing $c$. We repeat this procedure for $A'$ and a wall of $A'$, and so on.

We can again identify $Ht$ with the boundary of $U$ and then, by the description of $U$ and $\gamma_t$, we have
$$\gamma_t(Ht) = (q-1)^{-N}$$
where $N$ is the cardinal of the set of vertices with bad number.

Recall that by (5), one has,
$$Q - 1 = \frac{\log(q-1)}{\log a},$$
so we obtain
$$\gamma_t(Ht) = a^{-(Q-1)N} = (a^{-N})^{Q-1}.$$

We have to compare $a^{-N}$ and $\phi_l(t)$. For this, we consider the combinatorial Gromov product $\{.|.\}_x$ (see [1], page 256). Therefore, $\{\xi|t\}_x$ (respectively $\{t|\eta\}_x$, respectively $\{\xi|\eta\}_x$) is the number of walls of $A$ which intersects both $[x\xi)$ and $[xt)$ (respectively $[xt)$ and $[x\eta)$, respectively $[x\xi)$ and $[x\eta)$). Hence,
$$N = \{\xi|t\}_x + \{t|\eta\}_x - \{\xi|\eta\}_x.$$

Thus, by (5),
$$C^{-1}\frac{\delta_x(\xi,t)\delta_x(t,\eta)}{\delta_x(\xi,\eta)} \leq a^{-N} \leq C\frac{\delta_x(\xi,t)\delta_x(t,\eta)}{\delta_x(\xi,\eta)}.$$

But, $\dfrac{\delta_x(\xi,t)\delta_x(t,\eta)}{\delta_x(\xi,\eta)} = \dfrac{t(l-t)}{l}$ which is comparable to $\phi_l(t)$, and the Lemma 4 follows. $\square$

## 3. Proof of theorem 1

In this section we shall prove (2) for $\alpha = 1$.

The Poincaré inequality (2) is equivalent (see [7], Lemma 5.14) to the following inequality.

$$|u(\xi) - u(\eta)| \leq C\delta_x(\xi,\eta) \left(M_R\rho^\alpha(\xi) + M_R\rho^\alpha(\eta)\right)^{\frac{1}{\alpha}} \tag{10}$$

whenever $u : \partial I_{p,q} \longrightarrow \mathbb{R}^+$ is a continuous function in the ball $B_R$ of radius $R$, the points $\xi$ and $\eta$ are in $B_{C^{-1}R}$, $\rho$ is a very weak gradient of $u$ in $B_R$ and $M_R g$ is the maximal function defined by

$$M_R g(\xi) = \sup_{r<R} \frac{1}{\mu B(\xi,r)} \int_{B(\xi,r)} g\,d\mu.$$

We now begin the proof of (10) for $\alpha = 1$. The notations are the same as in the previous section.

By definition of the very weak gradient, for every $h \in H$,

$$|u(\xi) - u(\eta)| \leq \int_0^l \rho(ht)dt \text{ where } l = \delta_x(\xi,\eta). \tag{11}$$

Thus, integrating the last inequality, we obtain

$$\begin{aligned}
|u(\xi) - u(\eta)| &\leq \int_H \int_0^l \rho(ht)dt\,dh \text{ ( because $h$ is a probability)} \tag{12}\\
&= \int_0^l \left(\int_{Ht} \rho \Pi_{t*}(dh)\right) dt \text{ (by Fubini's theorem)}\\
&= \int_0^l \left(\frac{1}{\gamma_t(Ht)} \int_{Ht} \rho\,d\gamma_t\right) dt \text{ (by Lemma 4)}\\
&= \int_0^l f(t)dt \text{ (where } f(t) = \fint_{Ht} \rho\,d\gamma_t). \tag{13}
\end{aligned}$$

We shall need an elementary lemma of real analysis.

**Lemma 5.** *Let $Q \geq 1$.*
*Then, there exists a constant $C > 0$ such that*

$$\frac{1}{l}\int_0^l f(t)dt \leq C\left(\sup_{r\leq l} \frac{1}{r^Q}\int_0^r \phi_l^{Q-1}(t)f(t)dt + \sup_{r\leq l} \frac{1}{r^Q}\int_{l-r}^l \phi_l^{Q-1}(t)f(t)dt\right) \tag{14}$$

*whenever $l > 0$, $f$ is a positive Borel function on $[0,l]$ and where, for every $t \in [0,l]$, $\phi_l(t) = \inf(t, l-t)$.*

*Proof.* Note that it is sufficient to prove (14) for $l = 1$ (to see this, make the change of variables $u = \dfrac{t}{l}$ in the integrals and note that $\phi_l(t) = l\phi_1(\frac{t}{l})$).

We cut the integral of $f$ in the following way.

$$\begin{aligned}
\int_0^1 f(t)dt &= \sum_{n\geq 1} \int_{2^{-n-1}}^{2^{-n}} f(t)dt + \sum_{n\geq 1} \int_{1-2^{-n}}^{1-2^{-n-1}} f(t)dt\\
&\stackrel{\text{def}}{=} \sum_{n\geq 2} I_n + \sum_{n\geq 2} J_n. \tag{15}
\end{aligned}$$

But, for every $n \geq 1$,

$$(16) \quad I_n = \int_{2^{-n-1}}^{2^{-n}} f(t) dt = 2^{Q-1} 2^{-n} \left( \frac{1}{2^{-nQ}} \int_{2^{-n-1}}^{2^{-n}} (2^{-n-1})^{Q-1} f(t) dt \right).$$

Note that $\phi_1(t) = t$ if $t \in [2^{-n-1}, 2^{-n}]$. Thus,

$$I_n \leq 2^{Q-1} 2^{-n} \left( \frac{1}{(2^{-n})^Q} \int_0^{2^{-n}} \phi_1^{Q-1}(t) f(t) dt \right)$$

$$(17) \quad \leq 2^{Q-1} 2^{-n} \left( \sup_{r \leq 1} \frac{1}{r^Q} \int_0^r \phi_1^{Q-1}(t) f(t) dt \right).$$

By the same argument, we have

$$(18) \quad J_n \leq 2^{Q-1} 2^{-n} \left( \sup_{r \leq 1} \frac{1}{r^Q} \int_{1-r}^1 \phi_1^{Q-1}(t) f(t) dt \right).$$

The inequality (14) follows from (15), (17) and (18). $\square$

We shall now finish the proof of (10).
Applying the inequality (14) to the second member of (13), we obtain

$$|u(\xi) - u(\eta)| \leq Cl \sup_{r \leq l} r^{-Q} \int_0^r \phi_l^{Q-1}(t) \gamma_t(Ht)^{-1} \int_{Ht} \rho \, d\gamma_t \, dt$$

$$+ Cl \sup_{r \leq l} r^{-Q} \int_{1-r}^1 \phi_l^{Q-1}(t) \gamma_t(Ht)^{-1} \int_{Ht} \rho \, d\gamma_t \, dt.$$

Therefore using the Lemma 4 and the Ahlfors-regularity of $\mu$, we have

$$|u(\xi) - u(\eta)| \leq C \delta_x(\xi, \eta) \sup_{r \leq l} \frac{1}{\mu(B(\xi, r))} \int_0^r \int_{Ht} \rho \, d\gamma_t \, dt$$

$$+ C \delta_x(\xi, \eta) \sup_{r \leq l} \frac{1}{\mu(B(\eta, r))} \int_{1-r}^1 \int_{Ht} \rho \, d\gamma_t \, dt.$$

Fubini's theorem yields

$$|u(\xi) - u(\eta)| \leq C \delta_x(\xi, \eta) \sup_{r \leq \delta_x(\xi,\eta)} \frac{1}{\mu(B(\xi, r))} \int_{B(\xi,r) \cap \mathcal{C}} \rho \, d\mu$$

$$+ C \delta_x(\xi, \eta) \sup_{r \leq \delta_x(\xi,\eta)} \frac{1}{\mu(B(\eta, r))} \int_{B(\eta,r) \cap \mathcal{C}} \rho \, d\mu$$

where $\mathcal{C}$ is the "cone" defined in section 2.
The inequality (10) follows easily.

## ACKNOWLEDGMENTS

We would like to thank Thierry Coulhon and Pierre Pansu for their interest in this work. We whish to express our gratitude to Frédéric Paulin for the pictures and for valuable conversations.

Institut Elie Cartan, Département de mathématiques, Université de Nancy I, BP 239, 54506 Vandoeuvre les Nancy, France
*E-mail address*: bourdon@iecn.u-nancy.fr

MSRI, 1000 Centennial Drive, Berkeley, CA 94720-5070, USA, and Département de mathématiques, Université de Cergy-Pontoise, 2 avenue Adolphe Chauvin, 95302 Cergy-Pontoise cédex, France
*E-mail address*: pajot@u-cergy.fr